\documentclass[11pt]{article}
\usepackage{amsfonts}
\usepackage{amsmath}
\usepackage{graphicx}
\usepackage{url}

\newcommand{\eq}[1]{(\ref{eq:#1})}

\newtheorem{theorem}{Theorem}[section]

\newtheorem{lemma}[theorem]{Lemma}
\newtheorem{proposition}[theorem]{Proposition}

\numberwithin{equation}{section}

\newcommand{\qed}{\rule{7pt}{7pt}}
\newenvironment{proof}{\noindent{\bf Proof}\hspace*{1em}}{\hfill\qed\vspace{0.125in}}

\newcommand{\R}{\mathbb{R}}

\newcommand{\rank}{\operatorname{rank}}
\newcommand{\trace}{\operatorname{trace}}
\newcommand{\prob}{\mathbb{P}}

\newcommand{\minimize}{\mbox{minimize}}

\newcommand{\st}{\mbox{subject to}}

\newcommand{\cA}{\mathcal{A}}
\newcommand{\cG}{\mathfrak{G}}
\newcommand{\bfA}{\mathbf{A}}

\newcommand{\E}{\mathbb{E}}

\newcommand{\vvec}{\operatorname{vec}}

\title{\LARGE \bf
Necessary and Sufficient Conditions for Success of the Nuclear Norm
Heuristic for Rank Minimization}

\author{Benjamin Recht\thanks{Center for the Mathematics of Information,
California Institute of Technology, 1200 E California Blvd,
Pasadena, CA {\tt\small brecht@ist.caltech.edu}}, Weiyu
Xu\thanks{Electrical Engineering, California Institute of
Technology, 1200 E California Blvd, Pasadena, CA {\tt\small
weiyu@systems.caltech.edu}}, and Babak Hassibi\thanks{Electrical
Engineering, California Institute of Technology, 1200 E California
Blvd, Pasadena, CA {\tt\small bhassibi@systems.caltech.edu}} }

\begin{document}

\maketitle

\begin{abstract}
Minimizing the rank of a matrix subject to constraints is a
challenging problem that arises in many applications in control
theory, machine learning, and discrete geometry. This class of
optimization problems, known as rank minimization, is NP-HARD, and
for most practical problems there are no efficient algorithms that
yield exact solutions. A popular heuristic algorithm replaces the
rank function with the nuclear norm---equal to the sum of the
singular values---of the decision variable. In this paper, we
provide a necessary and sufficient condition that quantifies when
this heuristic successfully finds the minimum rank solution of a
linear constraint set.  We additionally provide a probability
distribution over instances of the affine rank minimization problem
such that instances sampled from this distribution satisfy our
conditions for success with overwhelming probability provided the
number of constraints is appropriately large.  Finally, we give
empirical evidence that these probabilistic bounds provide accurate
predictions of the heuristic's performance in non-asymptotic
scenarios.
\end{abstract}

\noindent{\footnotesize AMS (MOC) Subject Classification 90C25;
90C59; 15A52.}


\section{Introduction}

Optimization problems involving constraints on the rank of matrices
are pervasive in applications. In Control Theory, such problems
arise in the context of low-order controller
design~\cite{ElGhaoui93,Mesbahi97}, minimal realization
theory~\cite{Fazel01}, and model reduction~\cite{Beck98}.  In
Machine Learning, problems in inference with partial
information~\cite{Rennie05}, multi-task
learning~\cite{Argyriou08},and manifold learning~\cite{Weinberger06}
have been formulated as rank minimization problems.  Rank
minimization also plays a key role in the study of embeddings of
discrete metric spaces in Euclidean space~\cite{Linial95}. In
certain instances with special structure, rank minimization problems
can be solved via the singular value decomposition or can be reduced
to the solution of a linear system~\cite{Mesbahi97,PKrank}. In
general, however, minimizing the rank of a matrix subject to convex
constraints is NP-HARD. The best exact algorithms for this problem
involve quantifier elimination and such solution methods require at
least exponential time in the dimensions of the matrix variables.

A popular heuristic for solving rank minimization problems in the
controls community is the ``trace heuristic'' where one minimizes
the trace of a positive semidefinite decision variable instead of
the rank (see, e.g.,~\cite{Beck98,Mesbahi97}). A generalization of
this heuristic to non-symmetric matrices introduced by Fazel
in~\cite{FazelThesis} minimizes the \emph{nuclear norm}, or the sum
of the singular values of the matrix, over the constraint set.  When
the matrix variable is symmetric and positive semidefinite, this
heuristic is equivalent to the trace heuristic, as the trace of a
positive semidefinite matrix is equal to the sum of its singular
values. The nuclear norm is a convex function and can be optimized
efficiently via semidefinite programming.  Both the trace heuristic
and the nuclear norm generalization have been observed to produce
very low-rank solutions in practice, but, until very recently,
conditions where the heuristic succeeded were only available in
cases that could also be solved by elementary linear
algebra~\cite{PKrank}.

The first non-trivial sufficient conditions that guaranteed the
success of the nuclear norm heuristic were provided
in~\cite{Recht07}. Focusing on the special case where one seeks the
lowest rank matrix in an affine subspace, the authors provide a
``restricted isometry'' condition on the linear map defining the
affine subspace which guarantees the minimum nuclear norm solution
is the minimum rank solution. Moreover, they provide several
ensembles of affine constraints where this sufficient condition
holds with overwhelming probability.  Their work builds on seminal
developments in ``compressed sensing'' that determined conditions
for when minimizing the $\ell_1$ norm of a vector over an affine
space returns the sparsest vector in that space (see,
e.g.,~\cite{Candes05, CRT06, Wakin07}). There is a strong
parallelism between the sparse approximation and rank minimization
settings. The rank of a diagonal matrix is equal to the number of
non-zeros on the diagonal. Similarly, the sum of the singular values
of a diagonal matrix is equal to the $\ell_1$ norm of the diagonal.
Exploiting the parallels, the authors in~\cite{Recht07} were able to
extend much of the analysis developed for the $\ell_1$ heuristic to
provide guarantees for the nuclear norm heuristic.

Building on a different collection of developments in compressed
sensing~\cite{DT1,DT2,Stojnic08}, we present a \emph{necessary} and
sufficient condition for the solution of the nuclear norm heuristic
to coincide with the minimum rank solution in an affine space.  The
condition characterizes a particular property of the null-space of
the linear map which defines the affine space.  We show that when
the linear map defining the constraint set is generated by sampling
its entries independently from a Gaussian distribution, the
null-space characterization holds with overwhelming probability
provided the dimensions of the equality constraints are of
appropriate size.  We provide numerical experiments demonstrating
that even when matrix dimensions are small, the nuclear norm
heuristic does indeed always recover the minimum rank solution when
the number of constraints is sufficiently large.  Empirically, we
observe that our probabilistic bounds accurately predict when the
heuristic succeeds.


\subsection{Main Results}

Let $X$ be an $n_1 \times n_2$ matrix decision variable.  Without
loss of generality, we will assume throughout that $n_1\leq n_2$.
Let $\mathcal{A}:\mathbb{R}^{n_1\times n_2}\rightarrow \R^m$ be a
linear map, and let $b \in \R^m$. The main optimization problem
under study is
\begin{equation}
\begin{array}{ll}
\minimize & \rank(X)\\
\st & \cA(X)=b\,.
\end{array}
\label{eq:min-rank-prob}
\end{equation}

This problem is known to be NP-HARD and is also hard to
approximate~\cite{Meka08}. As mentioned above, a popular heuristic
for this problem replaces the rank function with the sum of the
singular values of the decision variable. Let $\sigma_{i}(X)$ denote
the $i$-th largest singular value of $X$ (equal to the square-root
of the $i$-th largest eigenvalue of $XX^*$). Recall that the rank of
$X$ is equal to the number of nonzero singular values. In the case
when the singular values are all equal to one, the sum of the
singular values is equal to the rank.  When the singular values are
less than or equal to one, the sum of the singular values is a
convex function that is strictly less than the rank.  This sum of
the singular values is a unitarily invariant matrix norm, called the
\emph{nuclear norm}, and is denoted
\[
\|X\|_* := \sum_{i=1}^{r} \sigma_i(X)\,.
\]
This norm is alternatively known by several other names including
the Schatten $1$-norm, the Ky Fan norm, and the trace class norm.

As described in the introduction, our main concern is when the
optimal solution of \eq{min-rank-prob} coincides with the optimal
solution of
\begin{equation}
\begin{array}{ll}
\minimize & \|X\|_*\\
\st & \cA(X)=b\,.
\end{array}
\label{eq:min-nucnorm-prob}
\end{equation}
This optimization is convex, and can be efficiently solved via a
variety of methods including semidefinite programming
(see~\cite{Recht07} for a survey).

Whenever $m<n_1n_2$, the null space of $\cA$, that is the set of $Y$
such that $\cA(Y)=0$, is not empty. Note that $X$ is an optimal
solution for \eq{min-nucnorm-prob} if and only if for every $Y$ in
the null-space of $\cA$
\begin{equation}\label{eq:weak-condition}
    \|X+Y\|_* \geq \|X\|_*\,.
\end{equation}
The following theorem generalizes this null-space criterion to a
critical property that guarantees when the nuclear norm heuristic
finds the minimum rank solution of $\cA(X)=b$ for all values of the
vector $b$. Our main result is the following

\begin{theorem}\label{thm:iff} Let $X_0$ be the optimal
solution of \eq{min-rank-prob} and assume that $X_0$ has rank
$r<n_1/2$. Then
\begin{enumerate}
\item If for every $Y$ in the null space of $\cA$ and for every decomposition
\[ Y = Y_1+Y_2, \]
where $Y_1$ has rank $r$ and $Y_2$ has rank greater than $r$, it
holds that
\[ \|Y_1\|_* < \|Y_2\|_*, \]
then $X_0$ is the unique minimizer of~\eq{min-nucnorm-prob}.
\item Conversely, if the condition of part 1 does not hold, then there
exists a vector $b\in\R^m$ such that the minimum rank solution of
$\cA(X)=b$ has rank at most $r$ and is not equal to the minimum
nuclear norm solution.
\end{enumerate}
\end{theorem}

This result is of interest for multiple reasons.  First,  as shown
in~\cite{RechtCDC08}, a variety of the rank minimization problems,
including those with inequality and semidefinite cone constraints,
can be reformulated in the form of~\eq{min-rank-prob}. Secondly, we
now present a family of random equality constraints under which the
nuclear norm heuristic succeeds with overwhelming probability. We
prove both of the following two theorems by showing that $\cA$ obeys
the null-space criteria of Equation \eq{weak-condition} and
Theorem~\ref{thm:iff} respectively with overwhelming
probability.

Note that for a linear map $\cA:\R^{n_1\times n_2}\rightarrow\R^m$,
we can always find an $m \times n_1n_2$ matrix $\bfA$ such that
\begin{equation}\label{eq:gauss-lin}
    \cA(X) = \bfA\vvec{X}\, .
\end{equation}
In the case where $\bfA$ has entries sampled independently from a
zero-mean, unit-variance Gaussian distribution, then the null space
characterization of theorem~\ref{thm:iff} holds with overwhelming
probability provided $m$ is large enough. For simplicity of notation
in the theorem statements, we consider the case of square matrices.
These results can be then translated into rectangular matrices by
padding with rows/columns of zeros to make the matrix square.  We
define the random ensemble of $d_1\times d_2$ matrices
$\cG(d_1,d_2)$ to be the Gaussian ensemble, with each entry sampled
i.i.d. from a Gaussian distribution with zero-mean and variance one.
We also denote $\cG(d,d)$ by $\cG(d)$.

The first result characterizes when a particular low-rank matrix can
be recovered from a random linear system via nuclear norm
minimization.
\begin{theorem}[Weak Bound]\label{thm:weak-bound}
Let $X_0$ be an $n\times n$ matrix of rank $r=\beta n$. Let $\cA:
 \R^{n\times n}\rightarrow \R^{\mu n^2}$ denote the random linear transformation
\begin{equation*}
    \mathcal{A}(X) = \bfA \vvec(X)\,,
\end{equation*}
where $\bfA$ is sampled from $\cG(\mu n^2,n^2)$.  Then whenever
\begin{equation}\label{eq:weak-bound}
    \mu \geq 1 - \frac{64}{9\pi^2}
\left((1-\beta)^{3/2}-\beta^{3/2}\right)^2
\end{equation}
there exists a numerical constant $c_w(\mu,\beta)>0$ such that with
probability exceeding $1-e^{ -c_w(\mu,\beta) n^2}$,
\begin{equation*}
    X_0 = \arg \min \{ \|Z\|_*~:~ \cA(Z) = \cA(X_0)\}\,.
\end{equation*}
In particular, if $\beta$ and $\mu$ satisfy \eq{weak-bound}, then
nuclear norm minimization will recover $X_0$ from a random set of
$\mu n^2$ constraints drawn from the Gaussian ensemble almost surely
as $n\rightarrow \infty$.
\end{theorem}

The second theorem characterizes when the nuclear norm heuristic
succeeds at recovering \emph{all} low rank matrices.

\begin{theorem}[Strong Bound]\label{thm:strong-bound}
Let $\cA$ be defined as in Theorem~\ref{thm:weak-bound}. Define the
two functions
\begin{align*}
    f(\beta,\epsilon) &= \frac{8}{3\pi} \frac{
(1-\beta)^{3/2} - \beta^{3/2}-4\epsilon }{1+4\epsilon}\\
    g(\beta,\epsilon) &=
    \sqrt{2\beta(2-\beta)}\log\left(\frac{3\pi}{2\epsilon}\right)\,.
\end{align*}
Then there exists a numerical constant $c_s(\mu,\beta)>0$ such that
with probability exceeding $1-e^{-c_s(\mu,\beta) n^2}$, for all
$n\times n$ matrices $X_0$ of rank $r \leq \beta n$
\begin{equation*}
    X_0 = \arg \min \{ \|Z\|_*~:~ \cA(Z) = \cA(X_0)\}
\end{equation*}
whenever
\begin{equation}\label{eq:strong-bound}
\begin{aligned}
\mu \geq 1 -
\sup_{\stackrel{\epsilon>0}{f(\beta,\epsilon)-g(\beta,\epsilon)>0}}
&\left(f(\beta,\epsilon)-g(\beta,\epsilon) \right)^2\,.
\end{aligned}
\end{equation}
In particular, if $\beta$ and $\mu$ satisfy \eq{weak-bound}, then
nuclear norm minimization will recover all rank $r$ matrices from a
random set of $\mu n^2$ constraints drawn from the Gaussian ensemble
almost surely as $n\rightarrow \infty$.
\end{theorem}

\begin{figure}
  \centering
   \includegraphics[width=8cm]{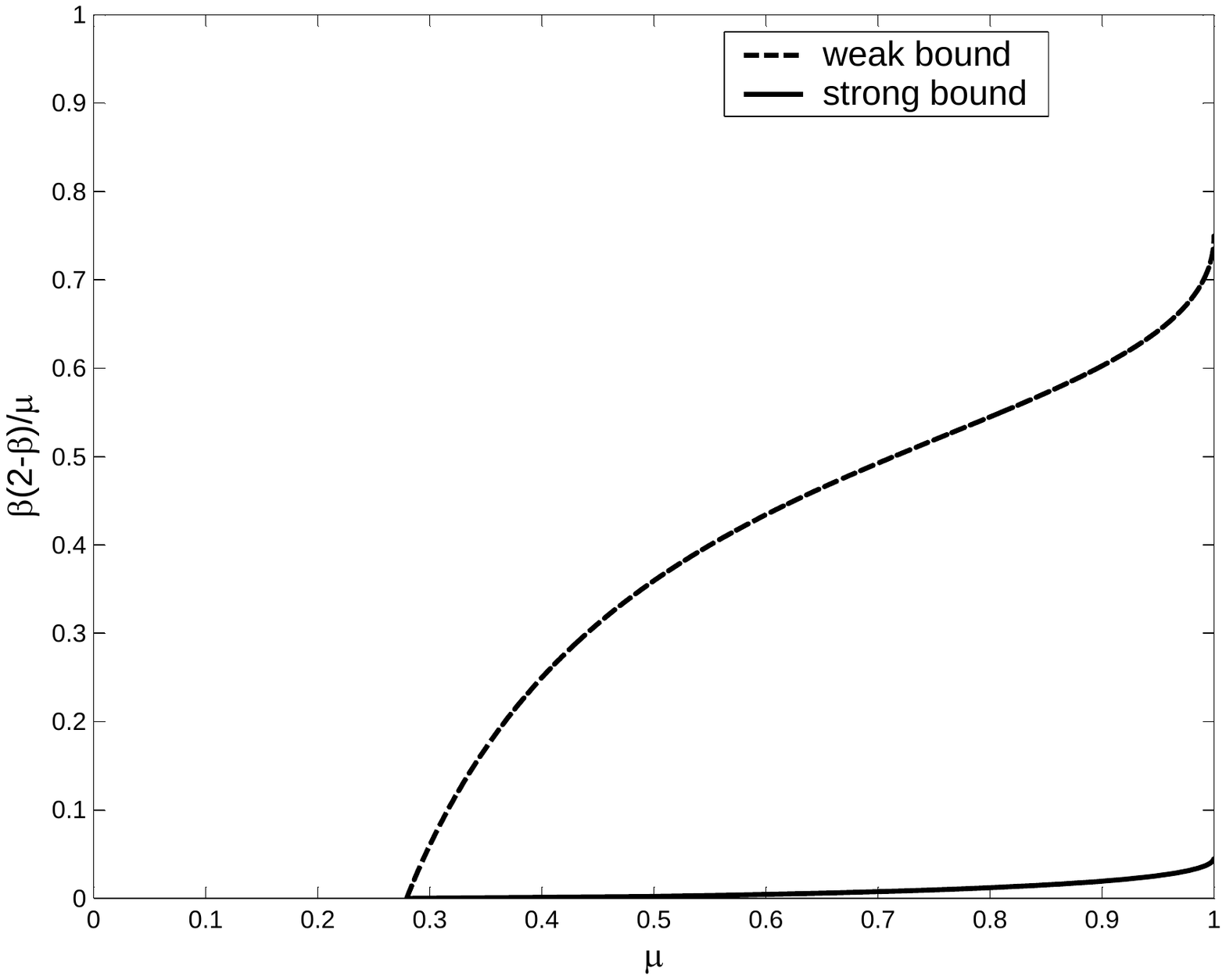}
   \caption{\small The Weak Bound~\eq{weak-bound} versus the Strong Bound~\eq{strong-bound}. } \label{fig:weak-v-strong}
\end{figure}

Figure~\ref{fig:weak-v-strong} plots the bound from
Theorems~\ref{thm:weak-bound} and~\ref{thm:strong-bound}. We call
\eq{weak-bound} the \emph{Weak Bound} because it is a condition that
depends on the optimal solution of~\eq{min-rank-prob}.  On the other
hand, we call \eq{strong-bound} the \emph{Strong Bound} as it
guarantees the nuclear norm heuristic succeeds no matter what the
optimal solution. The Weak Bound is the only bound that can be
tested experimentally, and, in Section~\ref{sec:experiments}, we
will show that it corresponds well to experimental data. Moreover,
the Weak Bound provides guaranteed recovery over a far larger region
of $(\beta,\mu)$ parameter space. Nonetheless, the mere existence of
a Strong Bound is surprising in of itself and results in a much
better bound than what was available from previous results
(c.f.,~\cite{Recht07}).

\subsection{Notation and Preliminaries}\label{sec:notation}

For a rectangular matrix $X \in \R^{n_1 \times n_2}$, $X^*$ denotes
the transpose of $X$.  $\vvec(X)$ denotes the vector in
$\R^{n_1n_2}$ with the columns of $X$ stacked on top of one and
other.

For vectors $v\in\R^d$, the only norm we will ever consider is the
Euclidean norm
\begin{equation*}
    \|v\|_{\ell_2} = \left(\sum_{i=1}^d v_i^2\right)^{1/2}\,.
\end{equation*}
On the other hand, we will consider a variety of matrix norms. For
matrices $X$ and $Y$ of the same dimensions, we define the inner
product in $\R^{n_1 \times n_2}$ as $\langle
X,Y\rangle:=\trace(X^*Y) = \sum_{i=1}^{n_1} \sum_{j=1}^{n_2} X_{ij}
Y_{ij}$. The norm associated with this inner product is called the
Frobenius (or Hilbert-Schmidt) norm $||\cdot||_F$. The Frobenius
norm is also equal to the Euclidean, or $\ell_2$, norm of the vector
of singular values, i.e.,
\[
\begin{split}
    \|X\|_F &:= \left(\sum_{i=1}^{r} {\sigma_i^2}\right)^\frac{1}{2}
    = \sqrt{\langle X,X\rangle}
    = \left( \sum_{i=1}^{n_1} \sum_{j=1}^{n_2} X_{ij}^2
    \right)^\frac{1}{2}
        \end{split}
\]
The operator norm (or induced 2-norm) of a matrix is equal to its
largest singular value (i.e., the $\ell_\infty$ norm of the singular
values):
\[
    \|X\| := \sigma_1(X).
\]
The nuclear norm of a matrix is equal to the sum of its singular
values, i.e.,
\[
\|X\|_* := \sum_{i=1}^{r} \sigma_i(X)\, \,.
\]
These three norms are related by the following inequalities which
hold for any matrix $X$ of rank at most $r$:
\begin{equation}
||X|| \leq ||X||_F \leq ||X||_* \leq \sqrt{r} ||X||_F \leq r ||X||.
\label{eq:ineqnorms}
\end{equation}

To any norm, we may associate a \emph{dual norm} via the following
variational definition
\begin{equation*}
    \|X\|_d = \sup_{\|Y\|_p=1} \langle Y, X\rangle\,.
\end{equation*}
One can readily check that the dual norm the Frobenius norm is the
Frobenius norm.  Less trivially, one can show that the dual norm of
the operator norm is the nuclear norm (See, for
example,~\cite{Recht07}).  We will leverage the duality between the
operator and nuclear norm several times in our analysis.


\section{Necessary and Sufficient Conditions}

We first prove our necessary and sufficient condition for success of
the nuclear norm heuristic.  We will need the following two
technical lemmas. The first is an easily verified fact.
\begin{lemma}\label{lemma:additive} Suppose $X$ and $Y$ are $n_1\times n_2$ matrices such that
$X^*Y=0$ and $XY^*=0$.  Then $\|X+Y\|_* = \|X\|_* +
\|Y\|_*$.\end{lemma}
Indeed, if $X^*Y=0$ and $XY^*=0$, we can find a
coordinate system in which
\begin{equation*}
    X = \left\|\left[\begin{array}{cc} A & 0\\0 & 0\end{array}\right]\right\|_*
    ~~\mbox{and}~~
    Y = \left\|\left[\begin{array}{cc} 0 & 0\\0 & B\end{array}\right]\right\|_*
\end{equation*} from which the lemma trivially follows.
The next Lemma allows us to exploit Lemma~\ref{lemma:additive} in
our proof.
\begin{lemma}\label{lemma:decomp}
Let $X$ be an $n_1\times n_2$ matrix with rank $r<\tfrac{n_1}{2}$
and $Y$ be an arbitrary $n_1\times n_2$ matrix. Let $P_X^c$ and
$P_X^r$ be the matrices that project onto the column and row spaces
of $X$ respectively. Then if $P_X^cYP_X^r$ has full rank, $Y$ can be
decomposed as
\begin{equation*} Y = Y_1+Y_2, \end{equation*} where $Y_1$ has rank $r$, and
\begin{equation*} \|X+Y_2\|_* = \|X\|_*+\|Y_2\|_*. \end{equation*}
\end{lemma}

\begin{proof} Without loss of generality, we can
write $X$ as
\[ X = \left[\begin{array}{cc} X_{11} & 0 \\ 0 & 0 \end{array}\right] , \]
where $X_{11}$ is $r\times r$ and full rank. Accordingly, $Y$
becomes
\[ Y = \left[\begin{array}{cc} Y_{11} & Y_{12} \\ Y_{21} & Y_{22} \end{array}\right] ,\]
where $Y_{11}$ is full rank since $P_X^rYP_X^c$ is. The
decomposition is now clearly
\[ Y = \underbrace{\left[\begin{array}{cc} Y_{11} & Y_{12} \\ Y_{21} &
Y_{21}Y_{11}^{-1}Y_{12} \end{array}\right]}_{Y_1} +
\underbrace{\left[\begin{array}{cc} 0 & 0
  \\ 0 & Y_{22}-Y_{21}Y_{11}^{-1}Y_{12} \end{array}\right]}_{Y_2}.
  \]  That $Y_1$ has rank $r$ follows from the fact that  the
rank of a block matrix is equal to the rank of a diagonal block plus
the rank of its Schur complement (see, e.g.,~\cite[\S
2.2]{HornJohnsonBook2}).  That $\|X_1+Y_2\|_* = \|X_1\|_* +
\|Y_2\|*$ follows from Lemma~\ref{lemma:additive}.
\end{proof}

We can now provide a proof of Theorem~\ref{thm:iff}.

\begin{proof} We begin by proving the converse. Assume the
condition of part 1 is violated, i.e., there exists some $Y$, such
that $\cA(Y) = 0$, $Y = Y_1+Y_2$, $\rank(Y_2)>\rank (Y_1)=r$, yet
$\|Y_1\|_* > \|Y_2\|_*$. Now take $X_0 = Y_1$ and $b = \cA(X_0)$.
Clearly, $\cA(-Y_2) = b$ (since $Y$ is in the null space) and so we
have found a matrix of higher rank, but lower nuclear norm.

For the other direction, assume the condition of part 1 holds.  Now
use Lemma~\ref{lemma:decomp} with $X=X_0$ and $Y = X_*-X_0$.  That
is, let $P_X^c$ and $P_X^r$ be the matrices that project onto the
column and row spaces of $X_0$ respectively and assume that
$P_{X_0}^c(X_*-X_0)P_{X_0}^r$ has full rank. Write $X_*-X_0 =
Y_1+Y_2$ where $Y_1$ has rank $r$ and $\|X_0+Y_2\|_* =
\|X_0\|_*+\|Y_2\|_*$. Assume further that $Y_2$ has rank larger than
$r$ (recall $r<n/2$). We will consider the case where
$P_{X_0}^c(X_*-X_0)P_{X_0}^r$ does not have full rank and/or $Y_2$
has rank less than or equal to $r$ in the appendix. We now have:
\begin{eqnarray*}
\|X_*\|_* & = & \|X_0+X_*-X_0\|_* \\
& = & \|X_0+Y_1+Y_2\|_* \\
& \geq & \|X_0+Y_2\|_*-\|Y_1\|_* \\
& = & \|X_0\|_*+\|Y_2\|_*-\|Y_1\|_*~~~\mbox{by
Lemma~\ref{lemma:decomp}. }
\end{eqnarray*}
But $\cA(Y_1+Y_2)=0$, so $\|Y_2\|_*-\|Y_1\|_*$ non-negative and
therefore $\|X_*\|_* \geq \|X_0\|_*$. Since $X_*$ is the minimum
nuclear norm solution, implies that $X_0 = X_*$. \end{proof}

For the interested reader, the argument for the case where
$P_{X_0}^r(X_*-X_0)P_{X_0}^c$ does not have full rank or $Y_2$ has
rank less than or equal to $r$ can be found in the appendix.

\section{Proofs of the Probabilistic Bounds}

We now turn to the proofs of the probabilistic
bounds~\ref{eq:weak-bound} and~\ref{eq:strong-bound}.   We first
provide a sufficient condition which implies the necessary and
sufficient null-space conditions.  Then, noting that the null space
of $\cA$ is spanned by Gaussian vectors, we use bounds from
probability on Banach Spaces to show that the sufficient conditions
are met.  The will require the introduction of two useful auxiliary
functions whose actions on Gaussian processes are explored in
Section~\ref{sec:aux-funs}.


\subsection{Sufficient Condition for Null-space Characterizations}\label{sec:sufficient} The following
theorem gives us a new condition that implies our necessary and
sufficient condition.

\begin{theorem}\label{thm:sufficient}
Let $\cA$ be a linear map of $n \times n$ matrices into $\R^m$.
Suppose that for every $Y$ in the null-space of $\cA$ and any
projection operators onto $r$-dimensional subspaces $P$ and $Q$ that
\begin{equation}\label{eq:suf-con}
    \|(I-P) Y (I-Q)\|_* \geq \|P Y Q\|_*\,.
\end{equation}
Then for every matrix $Z$ with row and column spaces equal to the
range of $Q$ and $P$ respectively,
\[
    \|Z+Y\|_* \geq \|Z\|_*
\]
for all $Y$ in the null-space of $\cA$.  In particular,
if~\ref{eq:suf-con} holds for every pair of projection operators $P$
and $Q$, then for every $Y$ in the null space of $\cA$ and for every
decomposition $Y = Y_1 + Y_2$ where $Y_1$ has rank $r$ and $Y_2$ has
rank greater than $r$, it holds that
\[
    \|Y_1\|_* \leq \|Y_2\|_*\,.
\]
\end{theorem}

We will need the following lemma

\begin{lemma}\label{lemma:nucbound}
For any block partitioned matrix
\[
    X = \left[\begin{array}{cc} A & B \\ C & D\end{array} \right]
\]
we have $\|X\|_* \geq \|A\|_* + \|D\|_*$.
\end{lemma}
\begin{proof}
This lemma follows from the dual description of the nuclear norm:
\begin{equation}\label{eq:upper-nuc}
\|X\|_* =  \sup \left\{ \left\langle \left[\begin{array}{cc} Z_{11}
&Z_{12}
\\Z_{21} & Z_{22}\end{array} \right], \left[\begin{array}{cc} A & B \\ C
& D\end{array} \right] \right\rangle \quad \bigg|\quad
\left\|\left[\begin{array}{cc} Z_{11} & Z_{12}
\\Z_{21} & Z_{22}\end{array} \right]\right\|=1\right\}\,.
\end{equation}
and similarly
\begin{equation}\label{eq:lower-nuc}
\|A\|_* + \|D\|_* = \sup \left\{ \left\langle
\left[\begin{array}{cc} Z_{11} &0
\\0 & Z_{22}\end{array} \right], \left[\begin{array}{cc} A & B \\ C
& D\end{array} \right] \right\rangle \quad \bigg|\quad
\left\|\left[\begin{array}{cc} Z_{11} & 0
\\0 & Z_{22}\end{array} \right]\right\|=1\right\}\,.
\end{equation}
Since \eq{upper-nuc} is a supremum over a larger set that
\eq{lower-nuc}, the claim follows.
\end{proof}

Theorem~\ref{thm:sufficient} now trivially follows

\begin{proof}[of Theorem~\ref{thm:sufficient}]
Without loss of generality, we may choose coordinates such that $P$
and $Q$ both project onto the space spanned by first $r$ standard
basis vectors. Then we may partition $Y$ as
\begin{equation*}
    Y = \left[\begin{array}{cc} Y_{11} & Y_{12} \\ Y_{21} & Y_{22}
    \end{array}\right]
\end{equation*}
and write, using Lemma~\ref{lemma:nucbound},
\begin{equation*}
   \|Y-Z\|_* - \|Z\|_* =   \left\| \left[\begin{array}{cc} Y_{11} - Z & Y_{12} \\ Y_{21} & Y_{22}
    \end{array}\right] \right\|_*  - \|Z\|_* \geq  \|Y_{11} - Z\|_* +
    \|Y_{22}\|_*  - \|Z\|_*\geq  \|Y_{22}\|_* - \|Y_{11}\|_*
\end{equation*}
which is non-negative by assumption.  Note that if the theorem holds
for all projection operators $P$ and $Q$ whose range has dimension
$r$, then $\|Z+Y\|_* \geq \|Z\|_*$ for all matrices $Z$ of rank $r$
and hence the second part of the theorem follows.
\end{proof}


\subsection{Proof of the Weak Bound}

Now we can turn to the proof of Theorem~\ref{thm:weak-bound}.  The
key observation in proving this lemma is the following
characterization of the null-space of $\cA$ provided by Stojnic et
al~\cite{Stojnic08}
\begin{lemma}\label{lemma:gaussian-null-space}
The null space of $\cA$ is identically distributed to the span of
$n^2(1-\mu)$ matrices $G_i$ where each $G_i$ is sampled i.i.d. from
$\cG(n)$.
\end{lemma}
This is nothing more than a statement that the null-space of $\cA$
is a random subspace.  However, when we parameterize elements in
this subspace as linear combinations of Gaussian vectors, we can
leverage Comparison Theorems for Gaussian processes to yield our
bounds.

Let $M=n^2(1-\mu)$ and let $G_1,\ldots, G_M$ be i.i.d. samples from
$\cG(n)$.  Let $X_0$ be a matrix of rank $\beta n$. Let $P_{X_0}$
and $Q_{X_0}$ denote the projections onto the column and row spaces
of $X_0$ respectively.  By theorem~\ref{thm:sufficient} and
Lemma~\ref{lemma:gaussian-null-space}, we need to show that for all
$v \in \R^M$,
\begin{equation}\label{eq:critical-bound}
    \|(I-P_{X_0}) \left( \sum_{i=1}^M v_i G_i \right)(I-Q_{X_0})\|_* \geq \|P_{X_0} \left(\sum_{i=1}^M v_i
    G_i\right) Q_{X_0}\|_*\,.
\end{equation}
That is, $\sum_{i=1}^M v_i G_i$ is an arbitrary element of the null
space of $\cA$, and this equation restates the sufficient condition
provided by Theorem~\ref{thm:sufficient}. Now it is clear by
homogeneity that we can restrict our attention to those $v\in \R^M$
with norm $1$.  The following crucial lemma characterizes when the
expected value of this difference is nonnegative

\begin{lemma}\label{lemma:inf-expectation-bound}
Let and $r=\beta n$ and suppose $P$ and $Q$ are projection operators
onto $r$-dimensional subspaces of $\R^n$.  For $i=1,\ldots, M$ let
$G_i$ be sampled from $\cG(n)$.  Then
\begin{equation}\label{eq:diff-norm-expectation}
\begin{aligned}
    \E\left[\inf_{\|v\|_{\ell_2}=1} \|(I-P)   \left(\sum_{i=1}^M v_i G_i\right) \right. & \left. (I-Q)\|_* -  \|P \left( \sum_{i=1}^M v_i
    G_i\right) Q\|_*\right]\\
    &\geq \left(\frac{8}{3\pi}+o(1)\right)\left( (1-\beta)^{3/2} - \beta^{3/2} \right)n^{3/2} -
    \sqrt{Mn}\,.
\end{aligned}
\end{equation}
\end{lemma}

We will prove this Lemma and a similar inequality required for the
proof the Strong Bound in Section~\ref{sec:aux-funs} below. But we
now show how using this Lemma and a concentration of measure
argument, we prove Theorem~\ref{thm:weak-bound}.

First note, that  if we plug in $M=(1-\mu)n^2$ and divide the right
hand side by $n^{3/2}$, the right hand side of
\eq{diff-norm-expectation} is non-negative if~\eq{weak-bound} holds.
To bound the probability that\eq{critical-bound} is non-negative, we
employ a powerful concentration inequality for the Gaussian
distribution bounding deviations of smoothly varying functions from
their expected value.

To quantify what we mean by smoothly varying, recall that a function
$f$ is \emph{Lipshitz} with respect to the Euclidean norm if there
exists a constant $L$ such that $|f(x)-f(y)| \leq L\|x-y\|_{\ell_2}$
for all $x$ and $y$.  The smallest such constant $L$ is called the
\emph{Lipshitz constant} of the map $f$. If $f$ is Lipshitz, it
cannot vary too rapidly.  In particular, note that if $f$ is
differentiable and Lipshitz, then $L$ is a bound on the norm of the
gradient of $f$.  The following theorem states that the deviations
of a Lipshitz function applied to a Gaussian random variable have
Gaussian tails.
\begin{theorem}\label{thm:concentration-of-measure}
Let $x$ be a normally distributed random vector and let $f$ be a
function with Lipshitz constant $L$. Then
\[
    \prob[|f(x)-\E[f(x)]|\geq t] \leq
    2\exp\left(-\frac{t^2}{2L^2}\right)\,.
\]
\end{theorem}
See~\cite{LedouxTalagrandBook} for a proof of this theorem with
slightly weaker constants and several references for more
complicated proofs that give rise to this concentration inequality.
The following Lemma bounds the Lipshitz constant of interest

\begin{lemma}~\label{lemma:lipshitz-inf}
For $i=1,\ldots, M$, let $X_i \in \R^{n_1 x n_1}$ and $Y_i \in
\R^{n_2\times n_2}$.  Define the function
\begin{equation*}
    F_I(X_1,\ldots, X_M, Y_1,\ldots, Y_M) = \inf_{\|v\|_{\ell_2}=1} \| \sum_{i=1}^M v_i X_i \|_* - \| \sum_{i=1}^M v_i
    Y_i\|_*\,.
\end{equation*}
Then the Lipshitz constant of $F_I$ is at most $\sqrt{n_1+n_2}$.
\end{lemma}

The proof of this lemma is straightforward and can be found in the
appendix. Using Theorem~\ref{thm:concentration-of-measure} and
Lemmas~\ref{lemma:inf-expectation-bound}
and~\ref{lemma:lipshitz-inf}, we can now bound
\begin{equation}\label{eq:inf-con}
\begin{aligned}
\prob\left[ \inf_{\|v\|_{\ell_2}=1} \|(I-P_{X_0})  \right. &
\left.\left(\sum_{i=1}^M v_i G_i\right) (I-Q_{X_0}) \|_*  -
\|P_{X_0} \left( \sum_{i=1}^M v_i
    G_i\right) Q_{X_0}\|_* \leq tn^{3/2}\right]\\
    &\leq
    \exp\left( -\tfrac{1}{2}\left(\frac{8}{3\pi}\left( (1-\beta)^{3/2} -
\beta^{3/2} \right) - \sqrt{1-\mu}-t\right)^2n^2 + o(n^2)\right)\,.
\end{aligned}
\end{equation}
Setting $t=0$ completes the proof of Theorem~\ref{thm:weak-bound}.
We will use this concentration inequality with a non-zero $t$ to
prove the Strong Bound.


\subsection{Proof of the Strong Bound}

The proof of the Strong Bound is similar to that of the Weak Bound
except we prove that \eq{critical-bound} holds for \emph{all}
operators $P$ and $Q$ that project onto $r$-dimensional subspaces.
Our proof will require an $\epsilon$-net for the projection
operators---a set of points such that any projection operator is
within $\epsilon$ of some element in the set.  We will show that if
a slightly stronger bound that \eq{critical-bound} holds on the
$\epsilon$-net, then \eq{critical-bound} holds for all choices of
row and column spaces.

Let us first examine how \eq{critical-bound} changes when we perturb
$P$ and $Q$. Let $P$, $Q$, $P'$ and $Q'$ all be projection operators
onto $r$-dimensional subspaces. Let $W$ be some $n\times n$ matrix
and observe that
\begin{align*}
 &\|(I-P') W(I-Q')\|_*  - \|P'W Q'\|_* - (\|(I-P) W(I-Q)\|_* - \|P
W
Q\|_*)\\
\leq&  \|(I-P) W(I-Q) - (I-P') W(I-Q')\|_* + \|P W Q - P' W Q'\|_*\\
\leq&   \|(I-P) W(I-Q) - (I-P') W(I-Q)\|_* + \|(I-P') W(I-Q) - (I-P') W(I-Q')\|_*\\
    &\qquad\qquad + \|P W Q - P' W  Q\|_* + \|P' W Q - P' W Q'\|_*\nonumber\\
\leq&   \|P-P'\|\|W\|_*\|I-Q\|  + \|I-P'\| \|W\|_*\|Q-Q'\| + \|P-P'\| \|W\|_* \|Q\| + \|P'\|\|W\|_*\|Q- Q'\|\nonumber\\
\leq&   2 (\|P-P'\|+\|Q-Q'\|) \|W \|_*\,.
\end{align*}
Here, the first and second lines follow from the triangle
inequality, the third line follows because $\|AB\|_* \leq
\|A\|\|B\|_*$, and the fourth line follows because $P$, $P'$, $Q$,
and $Q'$ are all projection operators.  Rearranging this inequality
gives
\[
\|(I-P') W(I-Q')\|_*  - \|P'W Q'\|_* \geq \|(I-P) W(I-Q)\|_* - \|P W
Q\|_*-   2 (\|P-P'\|+\|Q-Q'\|) \|W \|_*\,.
\]
As we have just discussed, if we can prove that with overwhelming
probability
\begin{equation}\label{eq:critical-bound-strong}
 \|(I-P) W(I-Q)\|_* - \|P W
Q\|_*-   4\epsilon \|W \|_* \geq 0
\end{equation}
for all $P$ and $Q$ in an $\epsilon$-net for the projection
operators onto $r$-dimensional subspaces, we will have proved the
Strong Bound.

To proceed, we need to know the size of an $\epsilon$-net.  The
following bound on such a net is due to Szarek.
\begin{theorem}[Szarek~\cite{SzarekHomog}]
Consider the space of all projection operators on $\R^n$ projecting
onto $r$ dimensional subspaces endowed with the metric
\[
    d(P,P') = \|P-P'\|
\]
Then there exists an $\epsilon$-net in this metric space with
cardinality at most
$\left(\tfrac{3\pi}{2\epsilon}\right)^{r(n-r/2-1/2)}$.
\end{theorem}

With this in hand, we now calculate the probability that for a given
$P$ and $Q$ in the $\epsilon$-net,
\begin{equation}
 \inf_{\|v\|_{\ell_2}=1} \|(I-P) \left(\sum_{i=1}^M v_i G_i\right)  (I-Q)\|_*- \|P \left(\sum_{i=1}^M v_i G_i\right) Q\|_*
 \geq 4\epsilon  \sup_{\|v\|_{\ell_2}=1} \left\|\sum_{i=1}^M v_i G_i
 \right\|_*\,.
\end{equation}

As we will show in Section~\ref{sec:aux-funs}, we can upper bound
the right hand side of this inequality using a similar bound as in
Lemma~\ref{lemma:inf-expectation-bound}.
\begin{lemma}\label{lemma:sup-expectation-bound}
For $i=1,\ldots, M$ let $G_i$ be sampled from $\cG(n)$.  Then
\begin{equation}\label{eq:diff-norm-expecation}
\begin{aligned}
    \E\left[\sup_{\|v\|_{\ell_2}=1} \|\sum_{i=1}^M v_i G_i\|_*\right] \leq \left(\frac{8}{3\pi}+o(1)\right)n^{3/2}
    + \sqrt{Mn}\,.
\end{aligned}
\end{equation}
\end{lemma}
Moreover, we prove the following in the appendix.
\begin{lemma}~\label{lemma:lipshitz-sup}
For $i=1,\ldots, M$, let $X_i \in \R^{n \times n}$ and define the
function
\[
    F_S(X_1,\ldots, X_M) = \sup_{\|v\|_{\ell_2}=1} \| \sum_{i=1}^M v_i X_i
    \|_*\,.
\]
Then the Lipshitz constant of $F_S$ is at most $\sqrt{n}$.
\end{lemma}

Using Lemmas~\ref{lemma:sup-expectation-bound}
and~\ref{lemma:lipshitz-sup} combined with
Theorem~\ref{thm:concentration-of-measure}, we have that
\begin{equation}\label{eq:sup-con}
    \prob\left[4\epsilon\sup_{\|v\|_{\ell_2}=1} \|\sum_{i=1}^M v_i G_i\|_* \geq
    t
    n^{3/2}
    \right] \leq \exp\left(-\tfrac{1}{2}\left( \frac{8}{3\pi} -
    \sqrt{1-\mu} + o(1) - \frac{t}{4\epsilon}\right)^2n^2\right)\,,
\end{equation}
and if we set the exponents of~\eq{inf-con} and \eq{sup-con} equal
to each other and solve for $t$, we find after some algebra and the
union bound
\begin{equation*}
\begin{aligned}
&\prob\left[ \inf_{\|v\|_{\ell_2}=1} \|(I-P) \left(\sum_{i=1}^M v_i
G_i\right) (I-Q)\|_*  - \|P\left(\sum_{i=1}^M v_i G_i\right) Q\|_*
\geq 4\epsilon \sup_{\|v\|_{\ell_2}=1} \|\sum_{i=1}^M v_i G_i\|_*
    \right] \\
     \geq& \prob\left[ \inf_{\|v\|_{\ell_2}=1} \|(I-P) \left(\sum_{i=1}^M v_i
G_i\right) (I-Q)\|_*  - \|P\left(\sum_{i=1}^M v_i G_i\right) Q\|_*>
    tn^{3/2}> 4\epsilon \sup_{\|v\|_{\ell_2}=1} \|\sum_{i=1}^M v_i G_i\|_*\right]\\
    \geq& 1 - \prob\left[ \inf_{\|v\|_{\ell_2}=1} \|(I-P) \left(\sum_{i=1}^M v_i
G_i\right) (I-Q)\|_*  - \|P\left(\sum_{i=1}^M v_i G_i\right)
Q\|_* < tn^{3/2}\right]\\
&\qquad\qquad\qquad\qquad\qquad\qquad - \prob\left[4\epsilon \sup_{\|v\|_{\ell_2}=1} \|\sum_{i=1}^M v_i G_i\|_* >  tn^{3/2}\right]\\
    \geq&
    1 - 2\exp\left(-\tfrac{1}{2}
    \left(\frac{8}{3\pi}\frac{(1-\beta)^{3/2}-\beta^{3/2}-4\epsilon}{1+4\epsilon}-\sqrt{1-\mu}\right)^2n^2+
    o(n^2)\right)\,.
\end{aligned}
\end{equation*}

Now, let $\Omega$ be an $\epsilon$-net for the set of projection
operators discussed above. Again by the union bound, we have that
\begin{equation}
\begin{aligned}
&\prob\left[\forall P,Q\, \inf_{\|v\|_{\ell_2}=1} \|(I-P)
\left(\sum_{i=1}^M v_i G_i\right) (I-Q)\|_*  - \|P\left(\sum_{i=1}^M
v_i G_i\right) Q\|_* \geq 4\epsilon \sup_{\|v\|_{\ell_2}=1}
\left\|\sum_{i=1}^M v_i
G_i Q\right\|_* \right] \\
& \leq 1-2\exp\left(-\left\{\tfrac{1}{2}
\left(\frac{8}{3\pi}\frac{(1-\beta)^{3/2}-\beta^{3/2}-4\epsilon}{1+4\epsilon}-\sqrt{1-\mu}\right)^2
+ \beta(2-\beta)\log\left(\frac{3\pi}{2\epsilon}\right)\right\}n^2 +
o(n)^2 \right)\,.
\end{aligned}
\end{equation}
Finding the parameters $\mu$, $\beta$, and $\epsilon$ that make the
terms multiplying $n^2$ negative completes the proof of the Strong
Bound.


\subsection{Comparison Theorems for Gaussian Processes and the Proofs
of Lemmas~\ref{lemma:inf-expectation-bound} and
\ref{lemma:sup-expectation-bound}}\label{sec:aux-funs}

Both of the two following Comparison Theorems provide sufficient
conditions for when the expected supremum or infimum of one Gaussian
process is greater to that of another.  Elementary proofs of both of
these Theorems and several other Comparison Theorems can be found in
\S 3.3 of~\cite{LedouxTalagrandBook}.

\begin{theorem}[Slepian's Lemma~\cite{Slepian62}]
Let $X$ and $Y$ by Gaussian random variables in $\R^N$ such that
\[
    \begin{cases}
        \E[X_iX_j] \leq \E[Y_iY_j] & \mbox{for all}~i\neq j\\
        \E[X_i^2] = \E[Y_i]^2 & \mbox{for all}~i
    \end{cases}
\]
Then
\[
  \E[\max_i Y_i] \leq \E[\max_i X_i]\,.
\]
\end{theorem}

\begin{theorem}[Gordan~\cite{Gordan85,Gordan88}]
Let $X=(X_{ij})$ and $Y=(Y_{ij})$ be Gaussian random vectors in
$\R^{N_1 \times N_2}$ such that
\[
    \begin{cases}
    \E[X_{ij}X_{ik}] \leq \E[Y_{ij}Y_{ik}] & \mbox{for all}~i,j,k\\
    \E[X_{ij}X_{lk}] \geq \E[Y_{ij}Y_{lk}] & \mbox{for all}~i\neq l~\mbox{and}~j,k\\
    \E[X_{ij}^2] = \E[X_{ij}^2] & \mbox{for all}~j,k
    \end{cases}
\]
Then
\[
    \E[\min_i \max_j  Y_{ij}] \leq \E[\min_i \max_j X_{ij}]\,.
\]
\end{theorem}

The following two lemmas follow from applications of these
Comparison Theorems.  We prove them in more generality than
necessary for the current work because both Lemmas are interesting
in their own right. Let $\|\cdot\|_p$ be any norm on $D\times D$
matrices and let $\|\cdot\|_d$ be its associated dual norm (See
Section~\ref{sec:notation}). Let us define the quantity
$\sigma(\|G\|_p)$ as
\begin{equation}
    \sigma(\|G\|_p) = \sup_{\|Z\|_d = 1} \|Z\|_F\,,
\end{equation}
and note that by this definition, we have
\[
    \sigma(\|G\|_p) = \sup_{\|Z\|_d = 1} \E\left[ \langle G,Z \rangle^2 \right]^{1/2}
\]
motivating the notation.

This first Lemma is now a straightforward consequence of Slepian's
Lemma
\begin{lemma}\label{lemma:sup-sup}
Let $\Delta>0$ and let $g$ be a Gaussian random vector in $\R^M$.
Let $G,G_1,\ldots,G_M$ be sampled i.i.d. from $\cG(D)$.  Then
\[
\E\left[\sup_{\|v\|_{\ell_2}=1}\sup_{\|Y\|_d=1} \Delta \langle
g,v\rangle + \left\langle \sum_{i=1}^M v_i G_i, Y \right\rangle
\right] \leq \E[\|G\|_p] +  \sqrt{M(\Delta^2 +
\sigma(\|G\|_p)^2)}\,.
\]
\end{lemma}

\begin{proof}

We follow the strategy used prove Theorem 3.20
in~\cite{LedouxTalagrandBook}.  Let $G,G_1,\ldots, G_M$ be sampled
i.i.d. from $\cG(D)$ and $g\in\R^M$ be a Gaussian random vector and
let $\gamma$ be a zero-mean, unit-variance Gaussian random variable.
For $v\in\R^M$ and $Y\in \R^{D\times D}$ define
\begin{align*}
Q_L(v,Y) &=  \Delta\langle g, v\rangle + \left\langle \sum_{i=1}^M
v_i G_i,
Y \right\rangle  + \sigma(\|G\|_p) \gamma\\
Q_R(v,Y) &=  \langle G, Y \rangle + \sqrt{\Delta^2 +
\sigma(\|G\|_p)^2} \langle g,v\rangle\,.
\end{align*}
Now observe that for any unit vectors in $\R^M$ $v$, $\hat{v}$ and
any $D\times D$ matrices $Y$, $\hat{Y}$ with dual norm $1$
\begin{equation*}
\begin{aligned}
& \E[Q_L(v,Y)Q_L(\hat{v},\hat{Y})] -
\E[Q_R(v,Y)Q_R(\hat{v},\hat{Y})]\\
=& \Delta^2 \langle v, \hat{v} \rangle + \langle v, \hat{v} \rangle
\langle Y, \hat{Y} \rangle + \sigma(\|G\|_p)^2 -  \langle Y, \hat{Y}
\rangle - (\Delta^2 + \sigma(\|G\|_p)^2)\langle v,
\hat{v} \rangle   \\
=& (\sigma(\|G\|_p)^2 - \langle Y, \hat{Y} \rangle)(1-\langle
v,\hat{v} \rangle)\,.
\end{aligned}
\end{equation*}
The difference in expectation is thus equal to zero if $v=\hat{v}$
and is greater than or equal to zero if $v\neq \hat{v}$.  Hence, by
Slepian's Lemma and a compactness argument (see
Proposition~\ref{prop:compact} in the Appendix),
\begin{equation*}
    \E\left[\sup_{\|v\|_{\ell_2}=1}\sup_{\|Y\|=1} Q_L(v,Y)\right]  \leq
\E\left[ \sup_{\|v\|_{\ell_2}=1}\sup_{\|Y\|=1} Q_R(v,Y)\right]
\end{equation*}
which proves the Lemma.
\end{proof}

The following lemma can be proved in a similar fashion

\begin{lemma}\label{lemma:inf-sup}
Let $\|\cdot\|_p$ be a norm on $\R^{D_1\times D_1}$ with dual norm
$\|\cdot\|_d$ and let $\|\cdot\|_b$ be a norm on $\R^{D_2 \times
D_2}$. Let $g$ be a Gaussian random vector in $\R^M$. Let
$G_0,G_1,\ldots,G_M$ be sampled i.i.d. from $\cG(D_1)$ and
$G'_1,\ldots,G'_M$ be sampled i.i.d. from $\cG(D_2)$. Then
\begin{equation*}
\begin{aligned}
\E\left[\inf_{\|v\|_{\ell_2}=1}\inf_{\|Y\|_{b}=1} \sup_{\|Z\|_{d}=1}
\right. & \left. \left\langle \sum_{i=1}^M v_i G_i, Z\right\rangle +
\left\langle \sum_{i=1}^M v_i
    G'_i, Y\right\rangle\right]\\
&\geq \E\left[\|G_0\|_p\right]  -
\E\left[\sup_{\|v\|_{\ell_2}=1}\sup_{\|Y\|_{b}=1}  \sigma(\|G\|_p)
\langle g, v \rangle + \left\langle \sum_{i=1}^M v_i G'_i,
Y\right\rangle\right]\,.
\end{aligned}
\end{equation*}
\end{lemma}

\begin{proof}
Define the functionals
\begin{align*}
    P_L(v,Y,Z)&= \left\langle \sum_{i=1}^M v_i G_i, Z\right\rangle +
\left\langle \sum_{i=1}^M v_i
    G'_i, Y\right\rangle + \gamma \sigma(\|G_0\|_p)\\
    P_R(v,Y,Z)&= \left\langle G_0, Z\right\rangle + \sigma(\|G_0\|_p)\langle g,
    v\rangle +
\left\langle \sum_{i=1}^M v_i G'_i, Y\right\rangle\,.
\end{align*}
Let $v$ and $\hat{v}$ be unit vectors in $\R^M$, $Y$ and $\hat{Y}$
be $D_2\times D_2$ matrices with $\|Y\|_b = \|\hat{Y}\|_b=1$, and
$Z$ and $\hat{Z}$ be $D_1\times D_1$ matrices with $\|Z\|_d =
\|\hat{Z}\|_d=1$.  Then we have
\begin{equation*}
\begin{aligned}
    &\E[P_L(v,Y,Z)P_L(\hat{v},\hat{Y},\hat{Z})] -
    \E[P_R(v,Y,Z)P_L(\hat{v},\hat{Y},\hat{Z})]\\
    =& \langle v, \hat{v} \rangle \langle Z, \hat{Z}\rangle + \langle v, \hat{v} \rangle \langle Y, \hat{Y} \rangle
    + \sigma(\|G_0\|_p)^2 - \langle Z, \hat{Z}\rangle -
    \sigma(\|G_0\|_p)^2\langle v, \hat{v} \rangle - \langle v, \hat{v} \rangle \langle Y, \hat{Y}
    \rangle\\
    =& (\sigma(\|G_0\|_p)^2 - \langle Z, \hat{Z}\rangle)(1 - \langle
    v, \hat{v} \rangle)\,.
\end{aligned}
\end{equation*}
The difference in expectations is greater than or equal to zero and
equal to zero when $v=\hat{v}$ and $Y=\hat{Y}$.  Hence, by Gordan's
Lemma and a compactness argument,
\begin{equation*}
    \E\left[\inf_{\|v\|_{\ell_2}=1}\inf_{\|Y\|_b=1}\sup_{\|Z\|_d=1} Q_L(v,Y,Z)\right]  \geq
\E\left[ \inf_{\|v\|_{\ell_2}=1}\inf_{\|Y\|_b=1}\sup_{\|Z\|_d=1}
Q_R(v,Y,Z)\right]
\end{equation*}
completing the proof.
\end{proof}

Together with Lemmas~\ref{lemma:sup-sup} and~\ref{lemma:inf-sup}, we
can prove the Lemma~\ref{lemma:inf-expectation-bound}.

\begin{proof}[of Lemma~\ref{lemma:inf-expectation-bound}]
For $i=1,\ldots,M$, let $G_i \in \cG((1-\beta)n)$ and $G'_i \in
\cG(\beta n)$.  Then
\begin{equation*}
\begin{aligned}
&    \E\left[\inf_{\|v\|_{\ell_2}=1} \left\|\sum_{i=1}^M v_i
G_i\right\|_*
  - \left\|  \sum_{i=1}^M v_i
    G'_i\right\|_*\right]\\
&=     \E\left[\inf_{\|v\|_{\ell_2}=1}\inf_{\|Y\|=1} \sup_{\|Z\|=1}
\left\langle \sum_{i=1}^M v_i G_i, Z\right\rangle  + \left\langle
\sum_{i=1}^M v_i
    G'_i, Y\right\rangle\right]\\
&\geq     \E\left[ \|G_0\|_*\right]  -
\E\left[\sup_{\|v\|_{\ell_2}=1}\sup_{\|Y\|=1} \sigma(\|G\|_*)
\langle g, v \rangle + \left\langle \sum_{i=1}^M v_i
    G'_i, Y\right\rangle \right]\\
    &\geq  \E\left[ \|G_0\|_*\right]  - \E\left[\|G'_0\|_*\right] -
    \sqrt{M}\sqrt{\sigma(\|G\|_*)^2 +\sigma(\|G'\|_*)^2}
\end{aligned}
\end{equation*}
where the first inequality follows from Lemma~\ref{lemma:inf-sup},
and the second inequality follows from Lemma~\ref{lemma:sup-sup}.

Now we only need to plug in the expected values of the nuclear norm
and the quantity $\sigma(\|G\|_*)$. Let $G$ be sampled from
$\cG(D)$. Then
\begin{equation}\label{eq:nuc-norm-expecation}
\E\|G\|_* = D\mathbb{E}\sigma_i = \frac{8}{3\pi} D^{3/2} + q(D)
\end{equation}
where $q(D)/D^{3/2}=o(1)$. The constant in from of the $D^{3/2}$
comes from integrating $\sqrt{\lambda}$ against the
Mar\v{c}enko-Pastur distribution (see,
e.g.,~\cite{Marcenko67,Bai99}):
\begin{equation*}
\frac{1}{2\pi}\int_{0}^{4} \sqrt{4-t} \, \, dt =
\frac{8}{3\pi}\approx 0.85\,.
\end{equation*}
Secondly, a straightforward calculation reveals
\begin{equation*}
\sigma(\|G\|_*)  = \sup_{\|H\| \leq 1} \|G\|_F = \sqrt{D} \,.
\end{equation*}
Plugging these values in with the appropriate dimensions completes
the proof.
\end{proof}

\begin{proof}[of Lemma~\ref{lemma:sup-expectation-bound}]
This lemma immediately follows from applying
Lemma~\ref{lemma:sup-sup} with $\Delta=0$ and from the calculations
at the end of the proof above. It is also an immediate consequence
of Lemma 3.21 from~\cite{LedouxTalagrandBook}.
\end{proof}


\section{Numerical Experiments}
\label{sec:experiments}

We now show that these asymptotic estimates hold even for small
values of $n$. We conducted a series of experiments for a variety of
the matrix sizes $n$, ranks $r$, and numbers of measurements $m$. As
in the previous section, we let $\beta = \tfrac{r}{n}$ and $\mu =
\tfrac{m}{n^2}$.  For a fixed $n$, we constructed random recovery
scenarios for low-rank $n\times n$ matrices.  For each $n$, we
varied $\mu$ between $0$ and $1$ where the matrix is completely
determined. For a fixed $n$ and $\mu$, we generated all possible
ranks such that $\beta(2-\beta)\leq \mu$. This cutoff was chosen
because beyond that point there would be an infinite set of matrices
of rank $r$ satisfying the $m$ equations.

For each $(n,\mu,\beta)$ triple, we repeated the following procedure
$10$ times. A matrix of rank $r$ was generated by choosing two
random $n\times r$ factors $Y_L$ and $Y_R$ with i.i.d. random
entries and setting $Y_0=Y_LY_R^*$. A matrix $\bfA$ was sampled from
the Gaussian ensemble with $m$ rows and $n^2$ columns. Then the
nuclear norm minimization
\begin{equation*}
    \begin{array}{ll}
        \mbox{minimize} &\quad \|X\|_*\\
        \mbox{subject to} &\quad \bfA\vvec{X} = \bfA\vvec{Y_0}
    \end{array}
\end{equation*}
was solved using the freely available software SeDuMi~\cite{sedumi}
using the semidefinite programming formulation described
in~\cite{Recht07}. On a 2.0 GHz Laptop, each semidefinite program
could be solved in less than two minutes for $40 \times 40$
dimensional $X$. We declared $Y_0$ to be recovered if
\[
    \|X-Y_0\|_F/\|Y_0\|_F<10^{-3}\,.
\]

\begin{figure*}
  \centering
  \begin{tabular}{cc}
    \includegraphics[width=8cm]{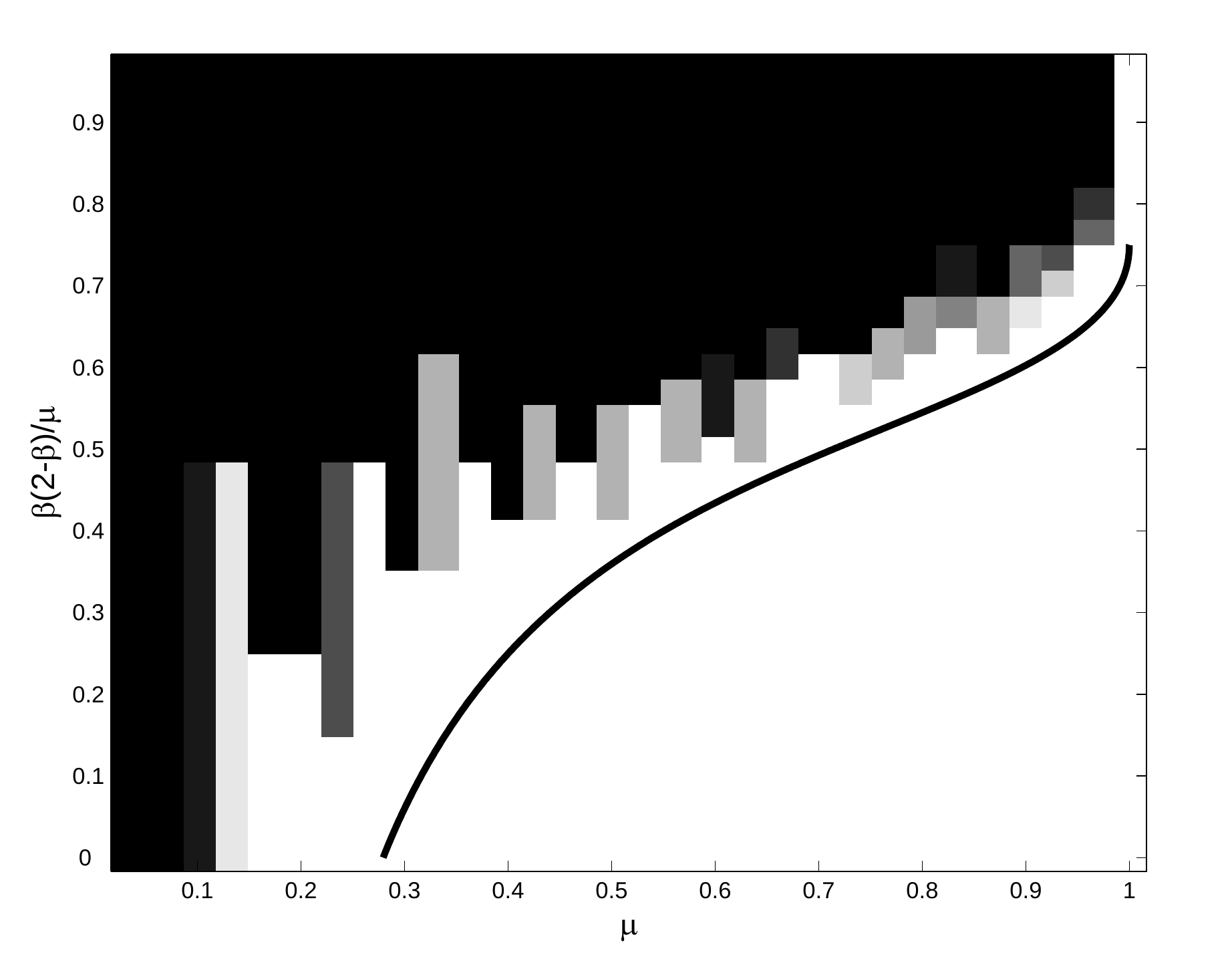}  &
    \includegraphics[width=8cm]{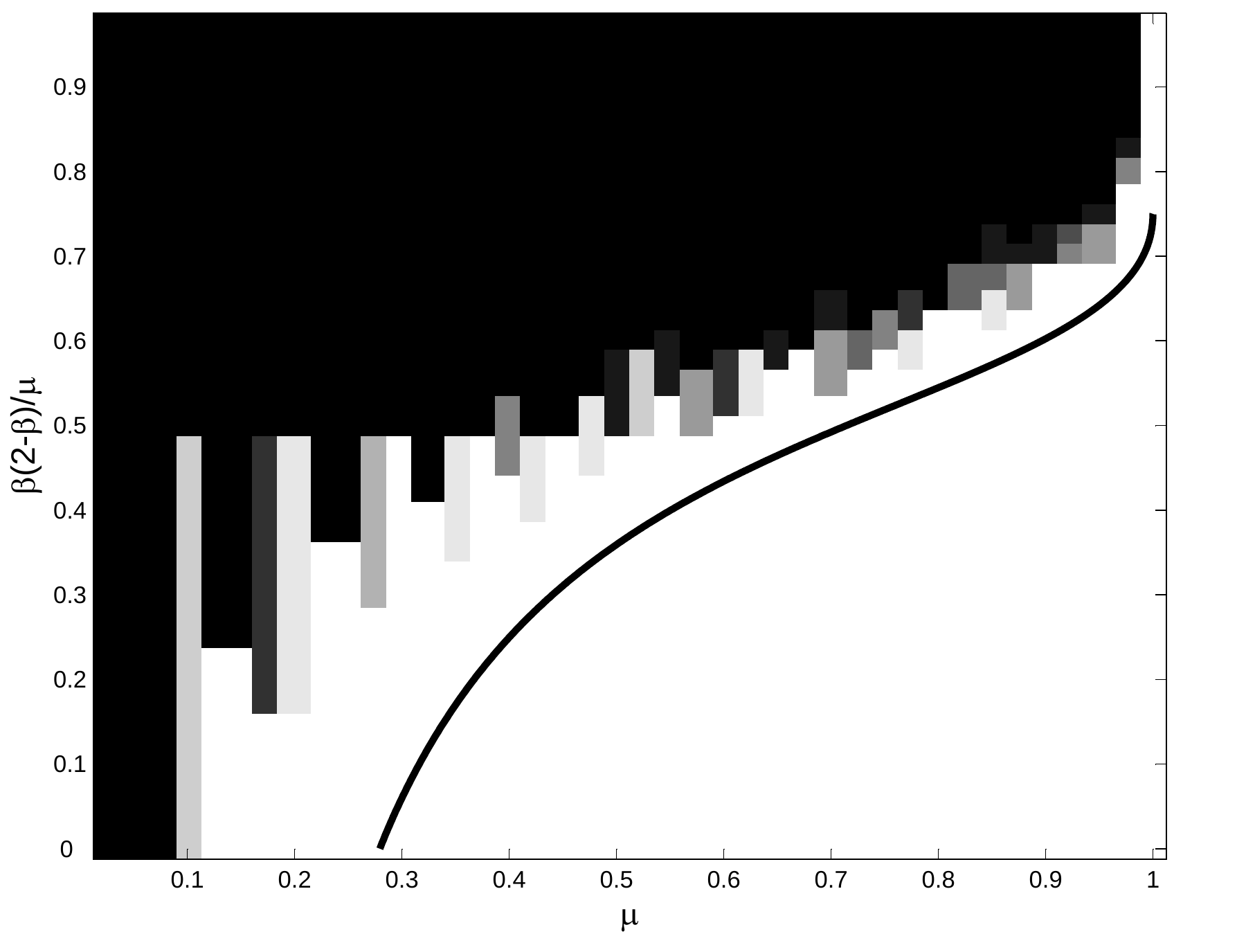}\\ (a) &
    (b)
      \end{tabular}
  \caption{\small
Random rank recovery experiments for  (a) $n=30$ and (b) $n=40$. The
color of each cell reflects the empirical recovery rate. White
denotes perfect recovery in all experiments, and black denotes
failure for all experiments.  In both frames, we plot the Weak
Bound~\eq{weak-bound}, showing that the predicted recovery regions
are contained within the empirical regions, and the boundary between
success and failure is well approximated for large values of
$\beta$.} \label{fig:phase-trans}
\end{figure*}

Figure~\ref{fig:phase-trans} displays the results of these
experiments for $n=30$ and $40$.  The color of the cell in the
figures reflects the empirical recovery rate of the $10$ runs
(scaled between $0$ and $1$).  White denotes perfect recovery in all
experiments, and black denotes failure for all experiments. It is
remarkable to note that not only are the plots very similar for
$n=30$ and $n=40$, but that the Weak Bound falls completely within
the white region and is an excellent approximation of the boundary
between success and failure for large $\beta$.

\if{0}
\section{Discussion and Conclusions}

We have presented a necessary and sufficient condition for the
nuclear norm heuristic (and hence also the trace heuristic) to find
the lowest rank solution of an affine set. Future work should also
investigate if the probabilistic analysis that provides the bounds
in Theorems~\ref{thm:weak-bound} and~\ref{thm:strong-bound} can be
further tightened at all. \fi

{\small
\bibliographystyle{IEEEtranS}
\bibliography{c:/latex/lowrank/manifesto/journal/brecht}
}

\appendix

\section{Appendix}

\subsection{Rank-deficient case of Theorem~\ref{thm:iff}}

As promised above, here is the completion of the proof of
Theorem~\ref{thm:iff}

\begin{proof}
In an appropriate basis, we may write
\[ X_0 = \left[\begin{array}{cc} X_{11} & 0 \\ 0 & 0 \end{array}\right] ~~~\mbox{and}~~~
X_*-X_0 = Y = \left[\begin{array}{cc} Y_{11} & Y_{12} \\ Y_{21} &
Y_{22}
\end{array}\right] \,.\] If $Y_{11}$ and $Y_{22}-Y_{21}Y_{11}^{-1}Y_{12}$ have full
rank, then all our previous arguments apply. Thus, assume that at
least one of them is not full rank. Nonetheless, it is always
possible to find an {\em arbitrarily small} $\epsilon >0$ such that
\[ Y_{11}+\epsilon I ~~~\mbox{and}~~~
\left[\begin{array}{cc} Y_{11}+\epsilon I & Y_{12} \\ Y_{21} &
Y_{22}+\epsilon I \end{array}\right]
\]
are full rank. This, of course, is equivalent to having
$Y_{22}+\epsilon I-Y_{21}(Y_{11}+\epsilon I)^{-1}Y_{12}$ full rank.
We can write
\begin{small}
\begin{align*}
\|X_*\|_* & =  \|X_0+X_*-X_0\|_* \\
& =  \left\|\left[\begin{array}{cc} X_{11} & 0 \\ 0 & 0
\end{array}\right] +\left[\begin{array}{cc} Y_{11} & Y_{12}
\\ Y_{21} & Y_{22} \end{array}\right] \right\|_*\\
& \geq  \left\| \left[\begin{array}{cc} X_{11}-\epsilon I & 0 \\ 0 &
Y_{22}-Y_{21}(Y_{11}+\epsilon I)^{-1}Y_{12} \end{array}\right]
\right\|_* - \left\| \left[\begin{array}{cc} Y_{11}+\epsilon I &
Y_{12}
\\ Y_{21} & Y_{21}(Y_{11}+\epsilon I)^{-1}Y_{12} \end{array}\right]
\right\|_* \\
& = \|X_{11}- \epsilon I\|_*  + \left\| \left[\begin{array}{cc} 0 &
0
\\ 0 & Y_{22}-Y_{21}(Y_{11}+\epsilon
I)^{-1}Y_{12} \end{array}\right] \right\|_* - \left\|
\left[\begin{array}{cc} Y_{11}+\epsilon I & Y_{12}
\\ Y_{21} & Y_{21}(Y_{11}+\epsilon I)^{-1}Y_{12} \end{array}\right]
\right\|_*\\
& \geq \|X_0\|_*  -r\epsilon + \left\| \left[\begin{array}{cc} \epsilon I-\epsilon I & 0 \\
0 & Y_{22}-Y_{21}(Y_{11}+\epsilon I)^{-1}Y_{12}
\end{array}\right] \right\|_* - \left\| \left[\begin{array}{cc} Y_{11}+\epsilon I & Y_{12}
\\ Y_{21} & Y_{21}(Y_{11}+\epsilon I)^{-1}Y_{12} \end{array}\right]
\right\|_*\\
& \geq  \|X_0\|_*-2r\epsilon \quad + \left\| \left[\begin{array}{cc}
-\epsilon I & 0 \\ 0 & Y_{22}-Y_{21}(Y_{11}+\epsilon I)^{-1}Y_{12}
\end{array}\right] \right\|_*  - \left\| \left[\begin{array}{cc} Y_{11}+\epsilon I & Y_{12}
\\ Y_{21} & Y_{21}(Y_{11}+\epsilon I)^{-1}Y_{12} \end{array}\right] \right\|_* \\
& \geq  \|X_0\|_*-2r\epsilon,
\end{align*}
\end{small}
where the last inequality follows from the condition of part 1 and
noting that
\begin{align*}
X_0-X_*&=\left[\begin{array}{cc} -\epsilon I & 0 \\
0 & Y_{22}-Y_{21}(Y_{11}+\epsilon I)^{-1}Y_{12} \end{array}\right] +
\left[\begin{array}{cc} Y_{11}+\epsilon I & Y_{12}
\\ Y_{21} &Y_{21}(Y_{11}+\epsilon I)^{-1}Y_{12}
\end{array}\right],
\end{align*}
lies in the null space of $\cA(\cdot)$ and the first matrix above
has rank more than $r$. But, since $\epsilon$ can be arbitrarily
small, this implies that $X_0 = X_*$. \end{proof}

\subsection{Lipshitz Constants of $F_I$ and $F_S$} We begin with the proof of
Lemma~\ref{lemma:lipshitz-sup} and then use this to estimate the
Lipshitz constant in Lemma~\ref{lemma:lipshitz-inf}.

\begin{proof}[of Lemma~\ref{lemma:lipshitz-sup}]
Note that the function $F_S$ is convex as we can write as a supremum
of a collection of convex functions
\begin{equation}\label{eq:FS-double-sup}
    F_S(X_1,\ldots,X_M) = \sup_{\|v\|_{\ell_2}=1} \sup_{\|Z\|<1} \langle \sum_{i=1}^M v_i  X_i, Z
    \rangle\,.
\end{equation}
The Lipshitz constant $L$ is bounded above by the maximal norm of a
subgradient of this convex function. That is, if we denote
$\bar{X}:=(X_1,\ldots,X_M)$, then we have
\begin{equation*}
    L \leq \sup_{\bar{X}} \sup_{\bar{Z} \in
    \partial F_S(\bar{X})} \left(\sum_{i=1}^M \|Z_i\|_F^2\right)^{1/2}\,.
\end{equation*}
Now, by \eq{FS-double-sup}, a subgradient of $F_S$ at $\bar{X}$ is
given of the form $(v_1 Z, v_2 Z,\ldots, v_M Z)$ where $v$ has norm
$1$ and $Z$ has operator norm $1$.  For any such subgradient
\begin{equation*}
    \sum_{i=1}^M \|v_i Z\|_F^2 = \|Z\|_F^2 \leq n
\end{equation*}
bounding the Lipshitz constant as desired.
\end{proof}

\begin{proof}[of Lemma~\ref{lemma:lipshitz-inf}]
For $i=1,\ldots,M$, let $X_i, \hat{X}_i\in \R^{n_1 x n_1}$, and
$Y_i,\hat{Y}_i \in \R^{n_2\times n_2}$. Let
\begin{equation*}
    w^* = \arg \min_{\|w\|_{\ell_2}=1}  \|\sum_{i=1}^M w_i \hat{X}_i\|_* -  \|\sum_{i=1}^M w_i
    \hat{Y}_i\|_*\,.
\end{equation*}
Then we have that
\begin{align*}
&    F_I(X_1,\ldots,X_M,Y_1,\ldots,Y_M) - F_I(\hat{X}_1,\ldots,\hat{X}_M,\hat{Y}_1,\ldots,\hat{Y}_M)\\
=& \left(\inf_{\|v\|_{\ell_2}=1} \|\sum_{i=1}^M v_i
    X_i\|_* -  \|\sum_{i=1}^M v_i Y_i\|_*\right) -
    \left(\inf_{\|w\|_{\ell_2}=1} \|\sum_{i=1}^M w_i
    \hat{X}_i\|_* -   \|\sum_{i=1}^M w_i \hat{Y}_i\|_*\right)\\
\leq & \|\sum_{i=1}^M w^*_i
    X_i\|_* -  \|\sum_{i=1}^M w^*_i Y_i\|_* -
    \|\sum_{i=1}^M w^*_i
    \hat{X}_i\|_* +   \|\sum_{i=1}^M w^*_i \hat{Y}_i\|_*\\
\leq & \|\sum_{i=1}^M w^*_i(X_i - \hat{X}_i)\|_* +  \|\sum_{i=1}^M
w^*_i (Y_i-\hat{Y}_i)\|_*\\
\leq & \sup_{\|w\|_{\ell_2}=1} \|\sum_{i=1}^M w_i(X_i -
\hat{X}_i)\|_* + \|\sum_{i=1}^M w_i (Y_i-\hat{Y}_i)\|_*  =
\sup_{\|w\|_{\ell_2}=1} \|\sum_{i=1}^M w_i \tilde{X}_i\|_* +
\|\sum_{i=1}^M w_i \tilde{Y}_i\|_*
\end{align*}
where $\tilde{X}_i = X_i - \hat{X}_i$ and $\tilde{Y}_i= Y_i -
\hat{Y}_i$.  This last expression is a convex function of
$\tilde{X}_i$ and $\tilde{Y}_i$,
\begin{equation*}
\sup_{\|w\|_{\ell_2}=1} \|\sum_{i=1}^M w_i \tilde{X}_i\|_* +
\|\sum_{i=1}^M w_i \tilde{Y}_i\|_* =
\sup_{\|w\|_{\ell_2}=1}\sup_{\|Z_X\|<1}\sup_{\|Z_Y\|<1} \langle
\sum_{i=1}^M w_i \tilde{X}_i, Z_X\rangle + \langle \sum_{i=1}^M w_i
\tilde{Y_i} Z_Y\rangle
\end{equation*}
with $Z_X$ $n_1\times n_2$ and $Z_Y$ $n_2 \times n_2$. Using an
identical argument as the one presented in the proof of
Lemma~\ref{lemma:lipshitz-sup}, we have that a subgradient of this
expression is of the form
\begin{equation*}
(w_1 Z_X, w_2 Z_X,\ldots, w_M Z_X, w_1 Z_Y, w_2 Z_Y,\ldots, w_M Z_Y)
\end{equation*}
where $w$ has norm $1$ and $Z_X$ and $Z_Y$ have operator norms $1$,
and thus
\begin{equation*}
     \sum_{i=1}^M \|w_iZ_X\|_F^2 + \|w_iZ_Y\|_F^2  =  \|Z_X\|_F^2 + \|Z_Y\|_F^2 \leq n_1+n_2
\end{equation*}
completing the proof.
\end{proof}

\subsection{Compactness Argument for Comparison Theorems}
\begin{proposition}\label{prop:compact}
Let $\Omega$ be a compact metric space with distance function
$\rho$. Suppose that $f$ and $g$ are real-valued function on
$\Omega$ such that $f$ is continuous and for any finite subset
$X\subset \Omega$
\begin{equation*}
    \max_{x\in X} f(x) \leq \max_{x\in X} g(x)\,.
\end{equation*}
Then
\begin{equation*}
    \sup_{x\in \Omega} f(x) \leq \sup_{x\in \Omega} g(x)\,.
\end{equation*}
\end{proposition}
\begin{proof}
Let $\epsilon>0$. Since $f$ is continuous and $\Omega$ is compact,
$f$ is uniformly continuous on $\Omega$.  That is, there exists a
$\delta>0$ such that for all $x,y\in\Omega$, $\rho(x,y)<\delta$
implies $|f(x)-f(y)|<\epsilon$. Let $X_\delta$ be a $\delta$-net for
$\Omega$.  Then, for any $x\in \Omega$, there is a $y$ in the
$\delta$-net with $\rho(x,y)<\delta$ and hence
\begin{equation*}
f(x) \leq f(y)+\epsilon \leq \sup_{z\in X_\delta} f(z)+\epsilon \leq
\sup_{z\in X_\delta} g(z) +\epsilon\leq \sup_{z\in \Omega} g(z) +
\epsilon\,.
\end{equation*}
Since this holds for all $x\in\Omega$ and $\epsilon>0$, this
completes the proof.
\end{proof}

\end{document}